# Kshirsagar–Tan independence property of beta matrices and related characterizations

KONSTANCJA BOBECKA and JACEK WESOŁOWSKI

*Wydział Matematyki i Nauk Informacyjnych, Politechnika Warszawska, Warszawa, Poland.
E-mail: bobecka@mini.pw.edu.pl; wesolo@mini.pw.edu.pl*

A new independence property of univariate beta distributions, related to the results of Kshirsagar and Tan for beta matrices, is presented. Conversely, a characterization of univariate beta laws through this independence property is proved. A related characterization of a family of $2 \times 2$ random matrices including beta matrices is also obtained. The main technical challenge was a problem involving the solution of a related functional equation.

## 1. Introduction

Characterizations of probability distributions by transformations preserving independence have been investigated by many authors. Celebrated theorems of Bernstein or Darmois and Skitovitch (for the normal law), of Lukacs (for the gamma law) and of Fisz (for the exponential law) are of this form. For beta laws, characterizations of this nature are typically related to neutrality properties of the Dirichlet distribution; see Fabius (1973), James and Mossiman (1980), Geiger and Heckerman (1997) and Bobecka and Wesołowski (2007). Matrix variate analogues of such characterizations are often much more demanding; see, for instance, characterizations of Wishart matrices due to Olkin and Rubin (1964), Casalis and Letac (1996), Geiger and Heckerman (2002), Bobecka and Wesołowski (2002) or characterizations of Wishart and matrix variate generalized inverse Gaussian matrices obtained in Letac and Wesołowski (2000) or Massam and Wesołowski (2006).

Our general goal is to develop characterizations related to independence properties of matrix variate beta distributions. A characterization of this type, related to matrix versions of neutrality properties, has recently been given by Hassairi and Regaig (2006). They extended to matrix random variables the result for univariate beta laws obtained in Seshadri and Wesołowski (2003). Here, we consider independence properties of beta matrices discovered by Kshirsagar (1961, 1972) and Tan (1969), which are essentially different than neutrality. This is carefully explained in Section 2. The study of the easiest possible case of $2 \times 2$ matrices led us to a new independence property for three univariate independent beta variables (for a recent review of the univariate beta distribution, one can consult Gupta and Nadarajah (2004)). The transformation preserving independence, we propose, seems to have no direct connection to neutrality (see the beginning of Section 3 and Section 5). The main results of this paper, given in Section 3, are characterizations: (1) of univariate beta laws through this new independence property; (2) of a family of $2 \times 2$ matrix variate distributions satisfying Kshirsagar–Tan independence properties. Somewhat unexpectedly, this latter family happens to be wider than the class of $2 \times 2$ matrix beta distributions. The proofs are based on a solution of a related functional equation. This equation is solved, under technical smoothness conditions, in Section 4. Concluding, in Section 5, we discuss similar characterizations for Wishart-type distributions and relations of our results to neutrality and to a representation of the Dirichlet distribution, as well as to uniqueness of solutions of some stochastic equations.

## 2. Independencies for beta matrices and beta distributions

Let $\mathcal{V}_n$ denote the Euclidean space of $n \times n$ symmetric matrices (with real entries) with the inner product defined as a trace of the product. The Lebesgue measure on $\mathcal{V}_n$ is fixed by assigning the unit mass to the unit cube in this space. Let $\mathcal{V}_n^+ \subset \mathcal{V}_n$ denote the cone of positive definite symmetric real $n \times n$ matrices. Additionally, denote by $\mathcal{M}_{n,m}$ the space of $n \times m$ real matrices. Let $e_n$ be the $n \times n$ identity matrix. Also, let $\mathcal{D}_n = \{x \in \mathcal{V}_n^+ : e_n - x \in \mathcal{V}_n^+\}$, an analogue of the interval $(0,1)$ in $\mathcal{V}_n^+$.

The matrix variate beta distribution $\beta_n(p,q)$ on $\mathcal{D}_n$ is defined by the density (with respect to the Lebesgue measure on $\mathcal{V}_n$) of the form

$$f(u) = \frac{(\det u)^{p-(n+1)/2}(\det(e-u))^{q-(n+1)/2}}{B_n(p,q)} I_{\mathcal{D}_n}(u),$$

where $p, q > \frac{n-1}{2}$ and $B_n(p,q)$ is the $n$-dimensional Euler beta function defined in terms of $n$-dimensional Euler gamma functions as

$$B_n(p,q) = \frac{\Gamma_n(p)\Gamma_n(q)}{\Gamma_n(p+q)}.$$

Recall that (see, e.g., Muirhead (1982), page 62)

$$\Gamma_n(p) = \pi^{n(n-1)/4} \prod_{i=1}^{n} \Gamma\left(p - \frac{i-1}{2}\right).$$

Let $X$ be a beta $n \times n$ random matrix $\beta_n(p,q)$. Let $T = [t_{ij}]$ be an upper triangular random matrix defined by $TT' = X$. Kshirsagar (1961, 1972) (see also Muirhead (1982),

Chapter 3.3 or Gupta and Nagar (2000), Chapter 5) proved that $t_{ii}^2$, $i=1,\ldots,n$, are independent beta variables, $t_{ii}^2 \sim B_I(p-\frac{i-1}{2},q)$, $i=1,\ldots,n$. In the course of the proof, he considered the block partitioning of $T$ with the dimensions of the blocks $1\times 1$, $1\times(n-1)$, $(n-1)\times 1$ and $(n-1)\times(n-1)$, as

$$T = \begin{bmatrix} t_{11} & \mathbf{t}' \\ 0 & T_{22} \end{bmatrix}$$

and showed that

$$t_{11}, T_{22} \quad \text{and} \quad \mathbf{v} = \frac{1}{\sqrt{1-t_{11}^2}}(e_{n-1}-T_{22}'T_{22})^{-1/2}\mathbf{t}$$

are independent. Moreover, it easily follows from his proof that $T_{22}'T_{22}$ has the matrix beta distribution $\beta_{n-1}(p-\frac{1}{2},q)$ and $\mathbf{v}$ has the density $f_{\mathbf{v}}(\mathbf{x}) \propto (1-\mathbf{x}'\mathbf{x})^{q-(n+1)/2}I_{S_{n-2}}(\mathbf{x})$, where $S_{n-2}$ is a unit ball in $\mathbb{R}^{n-1}$. This result was extended in Tan (1969), who considered the block partitioning of the matrix $X \sim \beta_n(p,q)$ with the dimensions of the blocks $r\times r$, $r\times s$, $s\times r$, $s\times s$,

$$X = \begin{bmatrix} X_{11} & X_{12} \\ X_{21} & X_{22} \end{bmatrix},$$

and proved that

$$X_{11}, \ X_{2\cdot 1} = X_{22} - X_{21}X_{11}^{-1}X_{12}, \quad V = (e_s - X_{2\cdot 1})^{-1/2}X_{21}(X_{11}^{-1}(e_r - X_{11})^{-1})^{1/2}$$

are independent random matrices. Moreover, it appears that $X_{11}$ is a beta matrix $\beta_r(p,q)$, $X_{2\cdot 1}$ is a beta matrix $\beta_s(p-\frac{r}{2},q)$ and $V$ has the density

$$f_V(x) \propto [\det(e_r - x'x)]^{q-(n+1)/2}I_{\mathcal{H}_{r,s}}(x), \qquad \text{where } \mathcal{H}_{r,s} = \{x \in \mathcal{M}_{s,r} : e_r - x'x \in \mathcal{V}_r^+\}.$$

Our aim is to study a converse problem, that is, we want to know if the independence conditions mentioned above characterize the matrix variate beta distribution. Since, in the general case, the problem seems to be very hard, here we study only the case of $2\times 2$ matrices. Let

$$X = \begin{bmatrix} X_{11} & X_{12} \\ X_{12} & X_{22} \end{bmatrix}$$

be a beta matrix $\beta_2(p,q)$. It then follows from the results of Kshirsagar (1961) and Tan (1969) that

$$X_{11}, X_{2\cdot 1} = X_{22} - \frac{X_{12}^2}{X_{11}}, V_1 = \frac{X_{12}}{\sqrt{(1-X_{11})(X_{11}-X_{11}X_{22}+X_{12}^2)}} \qquad \text{are independent,} \quad (2.1)$$

$$X_{22}, X_{1\cdot 2} = X_{11} - \frac{X_{12}^2}{X_{22}}, V_2 = \frac{X_{12}}{\sqrt{(1-X_{22})(X_{22}-X_{11}X_{22}+X_{12}^2)}} \qquad \text{are independent} \quad (2.2)$$

and both the triplets have the same distribution: $X_{11}$ and $X_{22}$ have the beta distribution $B_I(p,q)$, $X_{2\cdot1}$ and $X_{1\cdot2}$ have the beta distribution $B_I(p-\frac{1}{2},q)$, $V_1^2$ and $V_2^2$ have the beta distribution $B_I(\frac{1}{2},q-\frac{1}{2})$. It appears that the independence conditions (2.1) and (2.2) characterize a family of distributions which is wider than that of the beta matrix distributions.

First, instead of (2.1) and (2.2), we consider slightly weaker conditions, in the sense that the third element of both the triplets is squared. This allows us to rephrase the independence property completely in terms of univariate beta variables. As a matter of fact, an even more general property holds. This will be explained now.

Let

$$H = \{(x_1, x_2, x_3) \in (0,\infty)^3 : \min\{x_1 x_2, (1-x_1)(1-x_2)\} > x_3\}.$$

Denote by $B(p,q,r)$ the probability distribution on $H$ with the density

$$f(x_1, x_2, x_3) \propto (x_1 x_2 - x_3)^{p-1}[(1-x_1)(1-x_2) - x_3]^{q-1} x_3^{r-1} I_H(x_1, x_2, x_3),$$

where $p, q, r$ are positive numbers.

On $H$, define two mappings

$$\psi_i(x_1, x_2, x_3) = \left(x_i, \frac{x_1 x_2 - x_3}{x_i}, \frac{x_3}{(1-x_i)(x_i - x_1 x_2 + x_3)}\right), \qquad i = 1, 2. \quad (2.3)$$

Note that $\psi_i$, $i=1,2$, are bijections of $H$ onto $(0,1)^3$. Moreover, it can be checked that $\Psi = \psi_2 \circ \psi_1^{-1} : (0,1)^3 \to (0,1)^3$ is an involution (i.e., $\Psi = \Psi^{-1}$) and for any $(y_1, y_2, y_3) \in (0,1)^3$,

$$\begin{aligned}\Psi(y_1, y_2, y_3) \\ = \bigg(y_2 + (1-y_1)(1-y_2)y_3, \\ \frac{y_1 y_2}{y_2 + (1-y_1)(1-y_2)y_3}, \frac{y_1 y_3}{[1-(1-y_1)y_3][y_2 + y_3(1-y_2)]}\bigg).\end{aligned} \quad (2.4)$$

Let $(X_1, X_2, X_3)$ be a random vector having the distribution $B(p,q,r)$. A standard calculation using the Jacobian shows that $\psi_i(X_1, X_2, X_3)$ is a random vector with independent beta components

$$\psi_i(X_1, X_2, X_3) \sim B_I(p+r, q+r) \otimes B_I(p, q+r) \otimes B_I(r, q), \qquad i = 1, 2.$$

This observation can be rewritten as an independence property of a transformation of invariance beta variables.

**Theorem 1.** *If* $(Y_1, Y_2, Y_3) \sim B_I(p+r, q+r) \otimes B_I(p, q+r) \otimes B_I(r, q)$, *then*

$$\left(Y_2 + (1-Y_1)(1-Y_2)Y_3,\right.$$

$$\left. \frac{Y_1 Y_2}{Y_2 + (1 - Y_1)(1 - Y_2)Y_3}, \frac{Y_1 Y_3}{[1 - (1 - Y_1)Y_3][Y_2 + Y_3(1 - Y_2)]} \right) \stackrel{d}{=} (Y_1, Y_2, Y_3).$$

## 3. Characterizations

First, we study a characterization of beta laws which is a converse of Theorem 1. As already mentioned in the Introduction, to the best of our knowledge, all characterizations of beta variables existing in the literature which are based on transformations preserving independence can be related to neutrality properties of the Dirichlet distribution. The characterization of beta variables we give below is based on the independence of the components of the random vector defined in Theorem 1. It seems to have no connection to analogous results based on neutralities. This issue will be explained more carefully in Section 5.

**Theorem 2.** *Let $\mathbf{Y} = (Y_1, Y_2, Y_3)$ be a $(0,1)^3$-valued random vector with independent components having strictly positive, continuously differentiable density on $(0,1)^3$. Let $\psi_1$ and $\psi_2$ be the transformations defined in (2.3) and $\Psi = \psi_2 \circ \psi_1^{-1}$. Let*

$$\mathbf{Z} = (Z_1, Z_2, Z_3) = \Psi(Y_1, Y_2, Y_3)$$

*also have independent components.*

*There then exist positive numbers $p, q, r$ such that*

$$\mathbf{Y} \stackrel{d}{=} \mathbf{Z} \sim B_I(p+r, q+r) \otimes B_I(p, q+r) \otimes B_I(r, q).$$

**Proof.** Note that $\mathbf{Y} = \psi_1(\mathbf{X})$ and $\mathbf{Z} = \psi_2(\mathbf{X})$, where $\mathbf{X}$ is a random vector assuming values in $H$. The Jacobians of $\psi_i$, $i = 1, 2$, are, respectively,

$$|J_i| = \frac{1}{(1 - x_i)(x_i - x_1 x_2 + x_3)}, \qquad i = 1, 2.$$

Thus, the density $f$ of $\mathbf{X}$ can be expressed as

$$\begin{aligned}
f(\mathbf{x}) &= \frac{1}{(1 - x_1)(x_1 - x_1 x_2 + x_3)} \\
&\quad \times f_{Y_1}(x_1) f_{Y_2}\left(\frac{x_1 x_2 - x_3}{x_1}\right) f_{Y_3}\left(\frac{x_3}{(1 - x_1)(x_1 - x_1 x_2 + x_3)}\right) \\
&= \frac{1}{(1 - x_2)(x_2 - x_1 x_2 + x_3)} \\
&\quad \times f_{Z_1}(x_2) f_{Z_2}\left(\frac{x_1 x_2 - x_3}{x_2}\right) f_{Z_3}\left(\frac{x_3}{(1 - x_2)(x_2 - x_1 x_2 + x_3)}\right).
\end{aligned} \qquad (3.1)$$

We multiply (3.1) by $x_3$ and define

$$g_1(x) = \ln(f_{Y_1}(x)), \qquad g_2(x) = \ln(f_{Y_2}(x)), \qquad g_3(x) = \ln(xf_{Y_3}(x)),$$
$$g_4(x) = \ln(f_{Z_1}(x)), \qquad g_5(x) = \ln(f_{Z_2}(x)), \qquad g_6(x) = \ln(xf_{Z_3}(x)). \tag{3.2}$$

We then obtain

$$g_1(y_1) + g_2(y_2) + g_3(y_3) = g_4(z_1) + g_5(z_2) + g_6(z_3), \tag{3.3}$$

where $\psi_1(\mathbf{x}) = \mathbf{y}$ and $\psi_2(\mathbf{x}) = \mathbf{z}$. Note that $\mathbf{y} \in (0,1)^3$ and that

$$z_1 = y_2 + (1-y_1)(1-y_2)y_3, \qquad z_2 = \frac{y_1 y_2}{y_2 + (1-y_1)(1-y_2)y_3},$$
$$z_3 = \frac{y_1 y_3}{(1-(1-y_1)y_3)(y_2 + y_3 - y_2 y_3)}.$$

Note that equation (3.3) is the one we solve in Proposition 1 (Section 4). Hence, we conclude that there exist constants $p, q, r$ such that

$$f_{Y_1}(x) = f_{Z_1}(x) = B_1 x^{p+r-1}(1-x)^{q+r-1},$$
$$f_{Y_2}(x) = f_{Z_2}(x) = B_2 x^{p-1}(1-x)^{q+r-1},$$
$$f_{Y_3}(x) = f_{Z_3}(x) = B_3 x^{r-1}(1-x)^{q-1},$$

where the $B_i$'s are some normalizing constants. The integrability condition implies that $p, q, r$ are positive. $\square$

An alternative formulation of the above theorem is as follows.

**Theorem 3.** *Let $\mathbf{X} = (X_1, X_2, X_3)$ be an $H$-valued random vector having strictly positive, continuously differentiable density. If $\psi_1(\mathbf{X})$ and $\psi_2(\mathbf{X})$ have independent components, then there exist positive numbers $p, q, r$ such that $\mathbf{X} \sim B(p, q, r)$.*

If, instead, we use original independencies which hold for the beta $2 \times 2$ matrices, that is, (2.1) and (2.2), the functional equation of Proposition 1 (Section 4) leads to a family of distributions on $\mathcal{D}_2$ which is wider than that of beta matrix distributions.

**Theorem 4.** *Let $X$ be a $\mathcal{D}_2$-valued random matrix having a density which is continuously differentiable and strictly positive on $\mathcal{D}_2$. Assume that (2.1) and (2.2) hold. There then exist constants $a, b, c > 0$ such that the density of $X$ is of the form*

$$f_X\left(x = \begin{bmatrix} x_{11} & x_{12} \\ x_{12} & x_{22} \end{bmatrix}\right) = \frac{(\det x)^{a-1}(\det(e-x))^{b-1}|x_{12}|^{2c-1} I_{\mathcal{D}_2}(x)}{B(a,b)B(a+b,c)B(a+c,b+c)}. \tag{3.4}$$

**Proof.** Using (2.1) and (2.2) similarly as in the proof above, we get for the density $f_X$ the following two representations:

$$\frac{1}{\sqrt{(1-x_{11})(x_{11}-x_{11}x_{22}+x_{12}^2)}}$$
$$\times f_1(x_{11})f_2\left(\frac{x_{11}x_{22}-x_{12}^2}{x_{11}}\right)f_3\left(\frac{x_{12}}{\sqrt{(1-x_{11})(x_{11}-x_{11}x_{22}+x_{12}^2)}}\right)$$
$$=\frac{1}{\sqrt{(1-x_{22})(x_{22}-x_{11}x_{22}+x_{12}^2)}}$$
$$\times f_4(x_{22})f_5\left(\frac{x_{11}x_{22}-x_{12}^2}{x_{22}}\right)f_6\left(\frac{x_{12}}{\sqrt{(1-x_{22})(x_{22}-x_{11}x_{22}+x_{12}^2)}}\right),$$

where $f_i$, $i = 1, 2, \ldots, 6$, are respectively densities of $X_{11}, X_{2\cdot 1}, V_1, X_{22}, X_{1\cdot 2}, V_2$. We now proceed similarly as in the proof of Theorem 2, separately in two cases: $x_{12} > 0$ and $x_{12} < 0$.

For $x_{12} > 0$, denote $\sqrt{u} = x_{12}$. The above equation then takes on the form

$$f_1(x_{11})f_2\left(\frac{x_{11}x_{22}-u}{x_{11}}\right)\widetilde{f_3}\left(\frac{u}{(1-x_{11})(x_{11}-x_{11}x_{22}+u)}\right)$$
$$= f_4(x_{22})f_5\left(\frac{x_{11}x_{22}-u}{x_{22}}\right)\widetilde{f_6}\left(\frac{u}{(1-x_{22})(x_{22}-x_{11}x_{22}+u)}\right),$$

where $\widetilde{f_i}(x) = \sqrt{x}f_i(\sqrt{x})$, $i = 3, 6$. Thus, from Proposition 1, it follows that there exist positive (due to integrability) constants $a, b, c$ such that for $x \in (0, 1)$,

$$f_1(x) = f_4(x) = B_1 x^{a+c-1}(1-x)^{b+c-1},$$
$$f_2(x) = f_5(x) = B_2 x^{a-1}(1-x)^{b+c-1},$$
$$f_3(x) = f_6(x) = B_3 x^{2c-1}(1-x^2)^{b-1}.$$

For $x_{12} < 0$, denote $\sqrt{u} = -x_{12}$. The above equation then takes on the form

$$f_1(x_{11})f_2\left(\frac{x_{11}x_{22}-u}{x_{11}}\right)\hat{f_3}\left(\frac{u}{(1-x_{11})(x_{11}-x_{11}x_{22}+u)}\right)$$
$$= f_4(x_{22})f_5\left(\frac{x_{11}x_{22}-u}{x_{22}}\right)\hat{f_6}\left(\frac{u}{(1-x_{22})(x_{22}-x_{11}x_{22}+u)}\right),$$

where $\hat{f_i}(x) = \sqrt{x}f_i(-\sqrt{x})$, $i = 3, 6$. From Proposition 1, it follows that for $x \in (-1, 0)$,

$$f_3(x) = f_6(x) = B_3(-x)^{2c-1}(1-x^2)^{b-1}.$$

Thus, $f_3(x) = f_6(x) = B_3|x|^{2c-1}(1-x^2)^{b-1}I_{(-1,1)}(x)$.

Finally, the result follows from the representation of $f_X$ in terms of $f_1, f_2, f_3$. □

**Remark 1.** Note that the family of matrix variate distributions defined by (3.4) includes the beta matrix distribution: (3.4) for $c = \frac{1}{2}$ is the density of the beta matrix distribution $\beta_2(a + \frac{1}{2}, b + \frac{1}{2})$.

## 4. Functional equation

In this section, we solve the functional equation which was essential for proving the characterizations derived in Section 3.

**Proposition 1.** *Let $g_i$, $i = 1, 2, 3$, be continuously differentiable functions on $(0, 1)$ satisfying*

$$g_1(y_1) + g_2(y_2) + g_3(y_3)$$
$$= g_4(y_2 + (1 - y_1)(1 - y_2)y_3) + g_5\left(\frac{y_1 y_2}{y_2 + (1 - y_1)(1 - y_2)y_3}\right) \qquad (4.1)$$
$$+ g_6\left(\frac{y_1 y_3}{(1 - (1 - y_1)y_3)(y_2 + (1 - y_2)y_3)}\right)$$

*for any $y_1, y_2, y_3 \in (0, 1)$.*

*There then exist constants $\alpha, \beta, \gamma, A_i$, $i = 1, \ldots, 6$, such that $A_1 + A_2 + A_3 = A_4 + A_5 + A_6$ and*

$$\begin{aligned}
g_i(x) &= (\alpha + \gamma)\ln x + (\beta + \gamma)\ln(1 - x) + A_i, & i &= 1, 4, \\
g_i(x) &= \alpha \ln x + (\beta + \gamma)\ln(1 - x) + A_i, & i &= 2, 5, \\
g_i(x) &= \gamma \ln x + \beta \ln(1 - x) + A_i, & i &= 3, 6.
\end{aligned} \qquad (4.2)$$

**Proof.** Let $s = s(y_1, y_3) = 1 - (1 - y_1)y_3$, $t = t(y_2, y_3) = y_2 + (1 - y_2)y_3$. Taking derivatives of (4.1), once with respect to $y_1$, once with respect to $y_2$ and once with respect to $y_3$, yields the following three equations:

$$g_1'(y_1) = -(1 - y_2)y_3 g_4'(1 - (1 - y_2)s) \qquad (4.3)$$
$$+ \frac{y_2 t}{[1 - (1 - y_2)s]^2} g_5'\left(\frac{y_1 y_2}{1 - (1 - y_2)s}\right) + \frac{y_3(1 - y_3)}{s^2 t} g_6'\left(\frac{y_1 y_3}{st}\right),$$

$$g_2'(y_2) = s g_4'(1 - (1 - y_2)s) \qquad (4.4)$$
$$+ \frac{y_1(1 - s)}{[1 - (1 - y_2)s]^2} g_5'\left(\frac{y_1 y_2}{1 - (1 - y_2)s}\right) - \frac{y_1(1 - y_3)y_3}{st^2} g_6'\left(\frac{y_1 y_3}{st}\right),$$

$$g_3'(y_3) = (1 - y_1)(1 - y_2)g_4'(1 - (1 - y_2)s) \qquad (4.5)$$
$$- \frac{y_1(1 - y_1)y_2(1 - y_2)}{[1 - (1 - y_2)s]^2} g_5'\left(\frac{y_1 y_2}{1 - (1 - y_2)s}\right)$$

$$+ \frac{y_1[y_2 + (1-y_2)y_3(1-s)]}{s^2 t^2} g_6'\left(\frac{y_1 y_3}{st}\right).$$

Denote

$$G_j(x) = x g_j'(x), \qquad x \in (0,1), j = 1,2,3,4,5,6.$$

Then, upon multiplying (4.3), (4.4) and (4.5), respectively, by $y_1$, $y_2$ and $y_3$, we get

$$G_1(y_1) = -\frac{y_1(1-y_2)y_3}{1-(1-y_2)s} G_4(1-(1-y_2)s) \tag{4.6}$$

$$+ \frac{t}{1-(1-y_2)s} G_5\left(\frac{y_1 y_2}{1-(1-y_2)s}\right) + \frac{1-y_3}{s} G_6\left(\frac{y_1 y_3}{st}\right),$$

$$G_2(y_2) = \frac{y_2 s}{1-(1-y_2)s} G_4(1-(1-y_2)s) \tag{4.7}$$

$$+ \frac{(1-y_1)y_3}{1-(1-y_2)s} G_5\left(\frac{y_1 y_2}{1-(1-y_2)s}\right) - \frac{y_2(1-y_3)}{t} G_6\left(\frac{y_1 y_3}{st}\right)$$

and

$$G_3(y_3) = \frac{(1-y_2)(1-s)}{1-(1-y_2)s} G_4(1-(1-y_2)s)$$

$$- \frac{(1-y_2)(1-s)}{1-(1-y_2)s} G_5\left(\frac{y_1 y_2}{1-(1-y_2)s}\right) + \frac{y_2 + (1-y_2)y_3(1-s)}{st} G_6\left(\frac{y_1 y_3}{st}\right). \tag{4.8}$$

Letting $y_3 \to 0$ in (4.6), we conclude that the limit

$$\lim_{y_3 \to 0} G_6\left(\frac{y_1 y_3}{st}\right)$$

exists. We denote it by $G_6(0)$. (4.6) then yields

$$G_1(y_1) = G_5(y_1) + G_6(0). \tag{4.9}$$

Similarly, (4.6) for $y_3 \to 0$ gives

$$G_2(y_2) = G_4(y_2) - G_6(0). \tag{4.10}$$

Multiplying (4.7) by $1-y_2$ and adding the resulting equation to (4.8), we get

$$(1-y_2)G_2(y_2) + G_3(y_3) = (1-y_2)G_4(1-(1-y_2)s) + \frac{1-(1-y_2)s}{1-(1-y_1)y_3} G_6\left(\frac{y_1 y_3}{st}\right). \tag{4.11}$$

Letting $y_1 \to 1$ in (4.11), we obtain

$$(1-y_2)G_2(y_2) + G_3(y_3) = (1-y_2)G_4(y_2) + y_2 G_6\left(\frac{y_3}{t}\right).$$

Thus, by (4.10), we have

$$-(1-y_2)G_6(0) + G_3(y_3) = y_2 G_6\left(\frac{y_3}{t}\right),$$

and letting $y_2 \to 1$, we get

$$G_3(y_3) = G_6(y_3). \tag{4.12}$$

From (4.11), we have

$$G_6\left(\frac{y_1 y_3}{st}\right) = \frac{1-(1-y_1)y_3}{1-(1-y_2)s} \\ \times [(1-y_2)G_2(y_2) + G_3(y_3) - (1-y_2)G_4(1-(1-y_2)s)]. \tag{4.13}$$

We now plug (4.13) into (4.6) and (4.7). Multiplying the resulting equations by $(1-y_1)y_3$ and $-t$, respectively, and then adding them, we get

$$(1-y_1)y_3 G_1(y_1) - G_2(y_2) - (1-y_3)G_3(y_3) = -sG_4(1-(1-y_2)s). \tag{4.14}$$

Letting $y_2 \to 0$, we obtain

$$(1-y_1)y_3 G_1(y_1) - G_2(0) - (1-y_3)G_3(y_3) = -sG_4((1-y_1)y_3). \tag{4.15}$$

Note that the limit $\lim_{y\to 0} G_2(y)$ exists. We denote it by $G_2(0)$. Thus, (4.10) implies that the limit $\lim_{y\to 0} G_4(y)$ also exists. We denote it by $G_4(0)$. (4.10) then implies $G_4(0) = G_2(0) + G_6(0)$. Subtracting (4.15) from (4.14) and again using (4.10), we get

$$G_4(0) - G_4(y_2) - sG_4((1-y_1)y_3) = -sG_4(1-(1-y_2)s).$$

Let $x = (1-y_1)y_3 \in (0,1)$. Then, for any $x, y = y_2 \in (0,1)$ from the above equation, we obtain

$$G_4(0) - G_4(y) - (1-x)G_4(x) = -(1-x)G_4(y+x-yx). \tag{4.16}$$

Changing $x$ into $y$ and $y$ into $x$, we get

$$G_4(0) - G_4(x) - (1-y)G_4(y) = -(1-y)G_4(y+x-yx). \tag{4.17}$$

Multiplying (4.16) by $1-y$ and (4.17) by $1-x$ and then subtracting the resulting equations, we get

$$-yG_4(0) - (1-y_2)G_4(y) - (1-y)(1-x)G_4(x) \\ = -xG_4(0) - (1-x)G_4(x) - (1-x)(1-y)G_4(y).$$

Hence, for any $x, y \in (0,1)$,

$$-\frac{1}{y}G_4(0) + \frac{1-y}{y}G_4(y) = -\frac{1}{x}G_4(0) + \frac{1-x}{x}G_4(x)$$

and by separation of variables, we get

$$-\frac{1}{x}G_4(0) + \frac{1-x}{x}G_4(x) = C,$$

where $C$ is a constant. Thus, for any $x \in (0,1)$,

$$G_4(x) = -C + [C + G_4(0)]\frac{1}{1-x} \qquad (4.18)$$

and we get

$$g_4(x) = G_4(0)\ln x - (C + G_4(0))\ln(1-x) + A_4.$$

From (4.18) and (4.10), it follows that

$$G_2(x) = -(C + G_6(0)) + [C + G_4(0)]\frac{1}{1-x} \qquad (4.19)$$

and we obtain

$$g_2(x) = (G_4(0) - G_6(0))\ln x - (C + G_4(0))\ln(1-x) + A_2.$$

Plugging (4.18) and (4.19) into (4.15), we get

$$(1 - y_1)y_3 G_1(y_1) + G_6(0) - (1 - y_3)G_3(y_3) = -(1 - y_1)y_3 C.$$

Hence, for any $y_1, y_3 \in (0,1)$,

$$(1 - y_1)G_1(y_1) + (1 - y_1)C = -\frac{G_6(0)}{y_3} + \frac{1 - y_3}{y_3}G_3(y_3)$$

and again using separation of variables, we have

$$(1 - y)G_1(y) + (1 - y)C = D$$

and

$$-\frac{G_6(0)}{y} + \frac{1 - y}{y}G_3(y) = D,$$

where $D$ is a constant. Thus,

$$G_1(y) = -C + \frac{D}{1-y} \qquad (4.20)$$

and

$$G_3(y) = -D + [D + G_6(0)]\frac{1}{1-y}. \qquad (4.21)$$

Hence,

$$g_1(x) = (-C + D)\ln x - D\ln(1-x) + A_1$$

and

$$g_3(x) = G_6(0) \ln x - (D + G_6(0)) \ln(1-x) + A_3.$$

Using (4.9), we obtain

$$G_5(y) = -(C + G_6(0)) + \frac{D}{1-y} \qquad (4.22)$$

and using (4.12), we get

$$G_6(y) = -D + [D + G_6(0)] \frac{1}{1-y}. \qquad (4.23)$$

Hence,

$$g_5(x) = (-C + D - G_6(0)) \ln x - D \ln(1-x) + A_5$$

and

$$g_6(x) = G_6(0) \ln x - (D + G_6(0)) \ln(1-x) + A_6.$$

Inserting the formulas obtained for $g_i$, $i = 1, \ldots, 6$, into (4.1), we get that it must be

$$G_4(0) = D - C.$$

Thus, we get (4.2) where

$$\alpha = D - C - G_6(0), \qquad \beta = -D - G_6(0), \qquad \gamma = G_6(0). \qquad \square$$

## 5. Concluding remarks

The transformation $\Psi$ preserving a product beta distribution considered in this paper is defined in the three-variate situation. Describing its higher-dimensional version is an open problem.

On the other hand, it is natural to seek a multidimensional analogue of Theorem 4 by referring to matrix variate versions of independencies, as described in Section 2, for all or several block partitionings of the original matrix.

Recall that characterizations of Wishart matrices using independencies of $(X_{11}, X_{12})$ and $X_{2 \cdot 1}$ for all or several block partitionings of the $n \times n$ random matrix $X$ with $n > 2$ were proven in Geiger and Heckermann (2002) and in Massam and Wesołowski (2006). It is interesting to note that for $2 \times 2$ matrices, similarly as in Theorem 4 here, these independencies characterized wider families of random matrix distributions; see Geiger and Heckerman (1998) or Letac and Massam (2001).

Also, note that if $X$ is Wishart with a diagonal matrix parameter, then for any block partitioning of $X$, the following independence properties hold: $X_{11}$, $X_{21} X_{11}^{-1} X_{12}$ and $X_{2 \cdot 1}$ are independent and $X_{22}$, $X_{12} X_{22}^{-1} X_{21}$ and $X_{1 \cdot 2}$ are independent (see, e.g., Example 3.14 in Muirhead (1982)). Recently, a characterization using these independence properties for $2 \times 2$ matrices was proven in Seshadri and Wesołowski (2007). Analogously to our Theorem 4, it was shown there that the two independence conditions characterize a family of random matrices which is wider than that of Wishart matrices with diagonal matrix parameter. Again, the problem in higher dimensions remains open.

It is natural to compare the transformation $\Psi$ considered here to neutrality properties of the Dirichlet distribution (being consequences of the gamma variables representation), which define transformations preserving independence of beta variables. In particular, the one defined through the concept of complete neutrality (see Connor and Mosimann (1969), James and Mosimann (1980)) is worth analyzing. In the three-variate case which we are interested in here, the respective characterization has the following form. Let $\mathbf{Y} = (Y_1, Y_2, Y_3)$ be a $(0,1)^3$-valued random vector with independent components. Let

$$\mathbf{Z} = \left( \frac{Y_1}{1 - (1-Y_1)[Y_2 + (1-Y_2)Y_3]}, \frac{(1-Y_1)Y_2}{1 - (1-Y_1)(1-Y_2)Y_3}, (1-Y_1)(1-Y_2)Y_3 \right) \quad (5.1)$$

also have independent components. There then exist positive numbers $p, q, r, s$ such that

$$\mathbf{Y} \sim B_I(p, q+r+s) \otimes B_I(q, r+s) \otimes B_I(r, s)$$

and

$$\mathbf{Z} \sim B_I(p, s) \otimes B_I(q, p+s) \otimes B_I(r, p+q+s).$$

Let us emphasize that the transformation defined by (5.1) and $\Psi$ (see (2.4)) used in Theorems 1 and 2 are essentially different, though they may seem similar at first glance.

We conclude this section with two diverse and somewhat unexpected consequences of Theorem 1, one related to a representation of the bivariate Dirichlet distribution and the other related to stochastic equations and perpetuities.

*Remark 1.* Consider a random vector $(W_1, W_2)$ such that

$$\left( W_1, \frac{W_2}{1 - W_1} \right) \sim B_I(p, q+r) \otimes B_I(r, q).$$

Equivalently, $(W_1, W_2)$ has the Dirichlet distribution $D(p, r, q)$. By Theorem 1, looking at the last two coordinates of the random vector $\Psi(Y_1, Y_2, Y_3)$, we obtain the following representation:

$$(W_1, W_2) \stackrel{d}{=} U(V_1, V_2)$$

with

$$U = \frac{Y_1}{Y_2 + (1-Y_1)(1-Y_2)Y_3}, \quad V_1 = Y_2 \quad \text{and} \quad V_2 = \frac{(1-Y_1)Y_3}{1 - (1-Y_1)Y_3}.$$

It can easily be seen that $(V_1 + V_2)U \sim B_I(p+r, q)$, $V_1 \sim B_I(p, q+r)$ and $V_2 \sim B_{II}(r, p+q+r)$ are independent, where $B_{II}(r, p+q+r)$ denotes the second type of beta distribution

defined by the density

$$g(x) = \frac{x^{r-1}}{B(r, p+q+r)(1+x)^{p+q+2r}} I_{(0,\infty)}(x).$$

Also, note that $(\frac{V_1+V_2}{1+V_2}, UV_1, \frac{UV_2}{1-UV_1}) = \Psi(Y_1, Y_2, Y_3)$. Consequently, by Theorem 1, we conclude that $\frac{V_1+V_2}{1+V_2} \sim B_I(p+r, q+r)$ and $U(V_1, V_2) \sim D(p, r, q)$ are independent.

*Remark 2.* The results of Theorem 1 or 2 can be reviewed in the setting of stochastic equations. Namely, by comparing the respective coordinates of $\mathbf{Y}$ and $\Psi(\mathbf{Y})$, we observe that

**I.** The stochastic equation (for unknown $R$)

$$R \stackrel{d}{=} AR + B, \tag{5.2}$$

where $(A, B) \stackrel{d}{=} (-(1-Y_2)Y_3, Y_2 + (1-Y_2)Y_3)$ and $R$ are independent, has a solution $R \stackrel{d}{=} Y_1$;

**II.** The stochastic equation (for unknown $S$)

$$S \stackrel{d}{=} CS + D, \tag{5.3}$$

where $(C, D) \stackrel{d}{=} (\frac{(1-Y_1)Y_3}{Y_1}, \frac{1-(1-Y_1)Y_3}{Y_1})$ and $S$ are independent, has a solution $S \stackrel{d}{=} \frac{1}{Y_2}$;

**III.** The stochastic equation (for unknown $T$)

$$T \stackrel{d}{=} aT + b + \frac{c}{T}, \tag{5.4}$$

where $(a, b, c) \stackrel{d}{=} (\frac{Y_2}{Y_1}, \frac{1-2Y_2+Y_1Y_2}{Y_1}, -\frac{(1-Y_1)(1-Y_2)}{Y_1})$ and $T$ are independent, has a solution $T \stackrel{d}{=} \frac{1}{Y_3}$.

It follows from the theory of perpetuities (see, e.g., Vervaat (1979) and Goldie and Grübel (1996)) that for equations (5.2) and (5.3), these solutions are unique. However, it is not known if $T \stackrel{d}{=} \frac{1}{Y_3}$ is the unique solution of (5.4). Possibly, the reference most relevant to this uniqueness problem is the paper by Chamayou and Letac (1991), where many examples of stationary distributions for compositions of random functions are considered.

## Acknowledgements

This research was supported by the Ministry of Science and Higher Education of Poland (Project No. 1P03A022427).

# References

Bobecka, K. and Wesołowski, J. (2007). The Dirichlet distribution and process through neutralities. *J. Theor. Probab.* **20** 295–308. MR2324532

Bobecka, K. and Wesołowski, J. (2002). The Lukacs–Olkin–Rubin theorem without invariance of the "quotient". *Studia Mathematica* **152** 147–160. MR1916547

Casalis, M. and Letac, G. (1996). The Lukacs–Olkin–Rubin characterization of Wishart distributions on symmetric cones. *Ann. Statist.* **24** 763–786. MR1394987

Chamayou, J.F. and Letac, G. (1991). Explicit stationary distributions for compositions of random functions and products of random matrices. *J. Theor. Probab.* **4** 3–36. MR1088391

Connor, J.R. and Mosimann, J.E. (1969). Concepts of independence for proportions with a generalization of the Dirichlet distribution. *J. Amer. Statist. Assoc.* **64** 194–206. MR0240895

Fabius, J. (1973). Two characterizations of the Dirichlet distribution. *Ann. Statist.* **1** 583–587. MR0353531

Geiger, D. and Heckerman, D. (1997). A characterization of the Dirichlet distribution through global and local parameter independence. *Ann. Statist.* **25** 1344–1369. MR1447755

Geiger, D. and Heckerman, D. (1998). A characterization of the bivariate Wishart distribution. *Probab. Math. Statist.* **18** 119–131. MR1644057

Geiger, D. and Heckerman, D. (2002). Parameter priors for directed acyclic graphical models and the characterization of several probability distributions. *Ann. Statist.* **30** 1412–1440. MR1936324

Goldie, C.M. and Grübel, R. (1996). Perpetuities with thin tails. *Adv. in Appl. Probab.* **28** 463–480. MR1387886

Gupta, A.K. and Nadarajah, S. (2004). *Handbook of Beta Distribution and Its Applications*. New York: Dekker. MR2079703

Gupta, A.K. and Nagar, D.K. (2000). *Matrix Variate Distributions*. Boca Raton, FL: Chapman and Hall/CRC. MR1738933

Hassairi, A. and Regaig, O. (2006). Characterizations of the beta distribution on symmetric matrices. Preprint, 1-11.

James, I.R. and Mosimann, J.E. (1980). A new characterization of the Dirichlet distribution through neutrality. *Ann. Statist.* **8** 183–189. MR0557563

Kshirsagar, A.M. (1961). The non-central multivariate beta distribution. *Ann. Math. Statist.* **32** 104–111. MR0117833

Kshirsagar, A.M. (1972). *Multivariate Analysis*. New York: Dekker. MR0343478

Letac, G. and Massam, H. (2001). The normal quasi-Wishart distribution. In *Algebraic Methods in Statistics and Probability* (M.A.G. Viana and D.St.P. Richards, eds.). *AMS Contemporary Mathematics* **287** 231–239. MR1873678

Letac, G. and Wesołowski, J. (2000). An independence property for the product of GIG and gamma laws. *Ann. Probab.* **28** 1371–1383. MR1797878

Massam, H. and Wesołowski, J. (2006). The Matsumoto–Yor property and the structure of the Wishart distribution. *J. Multivariate Anal.* **97** 103–123. MR2208845

Muirhead, R.J. (1982). *Aspects of Multivariate Statistical Theory*. New York: Wiley. MR0652932

Tan, W.Y. (1969). Note on the multivariate and generalized multivariate beta distributions. *J. Amer. Statist. Assoc.* **64** 230–241. MR0240899

Olkin, I. and Rubin, H. (1964). Multivariate beta distributions and independence properties of the Wishart distribution. *Ann. Math. Statist.* **35** 261–269. MR0160297

Seshadri, V. and Wesołowski, J. (2003). Constancy of regressions for beta distributions. *Sankhyā A* **65** 284–291. MR2028900

Seshadri, V. and Wesołowski, J. (2007). More on connections between Wishart and matrix GIG distributions. *Metrika*. DOI: 10.1007/s00184-007-0154-3.

Vervaat, W. (1979). On a stochastic difference equation and a representation of non-negative infinitely divisible random variables. *Adv. in Appl. Probab.* **11** 750–783. MR0544194